\documentclass[12pt]{amsart}
\usepackage{amsmath, amsthm, amscd, amsfonts, amssymb, graphicx, color}

\makeatletter
\@namedef{subjclassname@2010}{%
  \textup{2010} Mathematics Subject Classification}
\makeatother

\newtheorem{theorem}{Theorem}[section]

\theoremstyle{definition}

\newtheorem{lemma}{Lemma}[section]

\newtheorem{remark}{Remark}[section]

\theoremstyle{theorem}
\newtheorem{other}{\bf Theorem}              

\newenvironment{pf}{\noindent{\emph{Proof.}}}{$\Box$ }
\newenvironment{Pf}{\noindent{\emph{Proof of}}}{$\Box$ }

\DeclareMathOperator{\defeq}{\overset{def}=}

\numberwithin{equation}{section}

\newcommand{\hol}{{\mathcal Hol}}
\DeclareMathOperator{\og}{O} 

\newcommand{\Aut}{{\rm Aut}}
\newcommand{\h}{\mathcal{H}}

\def\D{{\mathbb D}}

\begin{document}

\title{A Hankel matrix acting on spaces of analytic functions}

\author{Daniel Girela}
\address{An\'{a}lisis Matem\'{a}tico, Facultad de Ciencias, Universidad de M\'{a}laga, 29071 M\'{a}laga, Spain}
\email{girela@uma.es}
\author{Noel Merch\'{a}n}
\address{An\'{a}lisis Matem\'{a}tico, Facultad de Ciencias, Universidad de M\'{a}laga, 29071 M\'{a}laga, Spain}
\email{noel@uma.es}

\thanks{This research is supported in part by a grant from \lq\lq El Ministerio de
Econom\'{\i}a y Competitividad\rq\rq , Spain (MTM2014-52865-P) and
by a grant from la Junta de Andaluc\'{\i}a FQM-210. The second
author is also supported by a grant from \lq\lq El Ministerio de de
Educaci\'{o}n, Cultura y Deporte\rq\rq , Spain (FPU2013/01478).}



\begin{abstract} If $\mu $ is a positive Borel measure on the interval $[0, 1)$ we
let $\mathcal H_\mu $ be the Hankel matrix $\mathcal H_\mu =(\mu
_{n, k})_{n,k\ge 0}$ with entries $\mu _{n, k}=\mu _{n+k}$, where,
for $n\,=\,0, 1, 2, \dots $, $\mu_n$ denotes the moment of order $n$
of $\mu $. This matrix induces formally the operator
$$\mathcal{H}_\mu (f)(z)=
\sum_{n=0}^{\infty}\left(\sum_{k=0}^{\infty}
\mu_{n,k}{a_k}\right)z^n$$ on the space of all analytic functions
$f(z)=\sum_{k=0}^\infty a_kz^k$, in the unit disc $\D $. This is a
natural generalization of the classical Hilbert operator. In this
paper we improve the results obtained in some recent papers
concerning the action of the operators $H_\mu $ on Hardy spaces and
on M\"{o}bius invariant spaces.
\end{abstract}

\maketitle

\section{Introduction and main results}\label{intro}
We denote by $\mathbb D $ the unit disc in the complex plane
$\mathbb {C}$, and by $\hol (\mathbb D )$ the space of all analytic
functions in $\mathbb D$. We also let $H^p$ ($0<p\le \infty $) be
the classical Hardy spaces. We refer to \cite{D} for the notation
and results regarding Hardy spaces.
\par
\par
If\, $\mu $ is a finite positive Borel measure on $[0, 1)$ and $n\,
= 0, 1, 2, \dots $, we let $\mu_n$ denote the moment of order $n$ of
$\mu $, that is, $\mu _n=\int _{[0,1)}t^n\,d\mu (t),$ and we define
$\mathcal H_\mu $ to be the Hankel matrix $(\mu _{n,k})_{n,k\ge 0}$
with entries $\mu _{n,k}=\mu_{n+k}$. The matrix $\mathcal H_\mu $
can be viewed as an operator on spaces of analytic functions in the
following way: if\,
 $f(z)=\sum_{k=0}^\infty a_kz^k\in \hol (\D )$
we define
\begin{equation*}\label{H}
\mathcal{H}_\mu (f)(z)= \sum_{n=0}^{\infty}\left(\sum_{k=0}^{\infty}
\mu_{n,k}{a_k}\right)z^n,
\end{equation*}
whenever the right hand side makes sense and defines an analytic
function in $\D $.
\par\smallskip If $\mu $ is the Lebesgue measure on $[0,1)$ the matrix
$\mathcal H_\mu $ reduces to the classical Hilbert matrix \,
$\mathcal H= \left ({(n+k+1)^{-1}}\right )_{n,k\ge 0}$, which
induces the classical Hilbert operator $\h$ which has extensively
studied recently (see \cite{AlMonSa,Di,DiS,DJV,JK,LNP}). Other
related generalizations of the Hilbert operator have been considered
in \cite{GaGiPeSis} and \cite{PelRathg}.
\par The question of describing the measures $\mu $ for which the
operator $\mathcal H_\mu $ is well defined and bounded on distinct
spaces of analytic functions has been studied in a good number of
papers (see \cite{Bao-Wu, Ch-Gi-Pe, Ga-Pe2010, GM1, PV, Pow, Wi}).
Carleson measures play a basic role in these works.
\par If  $I\subset \partial\D$ is an
interval, $\vert I\vert $ will denote the length of $I$. The
\emph{Carleson square} $S(I)$ is defined as
$S(I)=\{re^{it}:\,e^{it}\in I,\quad 1-\frac{|I|}{2\pi }\le r <1\}$.
\par If $\, s>0$ and $\mu$ is a positive Borel  measure on  $\D$,
we shall say that $\mu $
 is an $s$-Carleson measure
  if there exists a positive constant $C$ such that
\[
\mu\left(S(I)\right )\le C{|I|^s}, \quad\hbox{for any interval
$I\subset\partial\D $}.
\]
\par  A $1$-Carleson
measure will be simply called a Carleson measure.
\par We recall that Carleson \cite{Carl}
proved that $H^p\,\subset \,L^p(d\mu )$ ($0<p<\infty $), if and only
if $\mu $ is a Carleson measure. This result was extended by Duren
\cite{Du:Ca} (see also \cite[Theorem\,\@9.\,\@4]{D}) who proved that
for $0<p\le q<\infty $, $H^p\subset L^q(d\mu )$ if and only if $\mu
$ is a $q/p$-Carleson measure.
\par
If $X$ is a subspace of $\hol (\D )$, $0<q<\infty $, and $\mu $ is a
positive Borel measure in $\mathbb D $, $\mu $ is said to be a
\lq\lq ${q}${\it {-Carleson measure for the space}} ${X}$\rq\rq \,
or an \lq\lq ${{(X, q)}}${\it-Carleson measure}\rq\rq \, if
$X\subset L^q(d\mu )$. The $q$-Carleson measures for the spaces
$H^p$, $0<p,q<\infty $ are completely characterized. The mentioned
results of Carleson and Duren can be stated saying that if $\,0<p\le
q<\infty $\, then a positive Borel measure $\mu $ in $\D$ is a
$q$-Carleson measure for $H^p$ if and only if $\mu$ is a
$q/p$-Carleson measure. Luecking \cite{Lu90} and Videnskii
\cite{Vid} solved the remaining case $0<q<p$. We mention
\cite{Bl-Ja} for a complete information on Carleson measures for
Hardy spaces.
\par
Following \cite{Zhao}, if $\mu$ is a positive Borel measure on $\D$,
$0\le \alpha <\infty $, and $0<s<\infty $ we say that $\mu$ is an
 $\alpha$-logarithmic $s$-Carleson measure if there exists a positive
 constant $C$ such that
 \[\frac{
\mu\left(S(I)\right )\left(\log \frac{2\pi }{\vert I\vert }\right
)^\alpha }{|I|^s}\le C, \quad\hbox{for any interval
$I\subset\partial\D $}.
\]
\par
A positive Borel measure $\mu $ on $[0, 1)$ can be seen as a Borel
measure on $\mathbb D$ by identifying it with the measure $\tilde
\mu $ defined by $$ \tilde \mu (A)\,=\,\mu \left (A\cap [0,1)\right
),\quad \text{for any Borel subset $A$ of $\mathbb D$}.$$  In this
way a positive Borel measure $\mu $ on $[0, 1)$ is an $s$-Carleson
measure if and only if there exists a positive constant $C$ such
that
\[
\mu\left([t,1)\right )\le C(1-t)^s, \quad 0\le t<1.
\]
We have a similar statement for $\alpha$-logarithmic $s$-Carleson
measures.
\par Widom \cite[Theorem\,\@3.\,\@1]{Wi} (see also
\cite[Theorem\,\@3]{Pow} and \cite[p.\,\@42,
Theorem\,\@7.\,\@2]{Pell}) proved that $\mathcal H_\mu $ is a
bounded operator from $H^2$ into itself if and only $\mu $ is a
Carleson measure. Galanopoulos and Pel\'{a}ez \cite{Ga-Pe2010}
studied the operators $\mathcal H_\mu $ acting on $H^1$ and
Chatzifountas, Girela and Pel\'{a}ez \cite{Ch-Gi-Pe} studied the
action of $\mathcal H_\mu $ on $H^p$, $0<p<\infty $. \par A key
ingredient in \cite{Ga-Pe2010} and \cite{Ch-Gi-Pe} is obtaining an
integral representation of $\mathcal H_\mu $. If $\mu $ is as above,
we shall write throughout the paper
\begin{equation*}I_\mu
(f)(z)=\int_{[0,1)}\frac{f(t)}{1-tz}\,d\mu (t),\end{equation*}
whenever the right hand side makes sense and defines an analytic
function in $\D $. It turns out that the operators $H_\mu $ and
$I_\mu$ are closely related. Indeed, some of the results obtained in
\cite{Ga-Pe2010} and \cite{Ch-Gi-Pe} are the following ones:
\begin{other}[\cite{Ga-Pe2010}]\label{ImuHmuH1} Let
$\mu $ be a positive Borel measure on $[0,1)$. Then:
\begin{itemize}
\item[(i)] The operator $I_\mu $ is well defined on $H^1$ if and only
if $\mu $ is a Carleson measure.
\item[(ii)] If $\mu $ is a Carleson measure, then the
operator $\mathcal H_\mu $ is also well defined on $H^1$ and,
furthermore,
$$\mathcal H_\mu (f)\,=\, I_\mu (f),\quad\text{for every $f\in
H^1$}.$$
\item[(iii)] The operator $I_\mu $ is a bounded operator from $H^1$
into itself if and only if $\mu $ is a $1$-logarithmic $1$-Carleson
measure.
\end{itemize}
\end{other}

\begin{other}[\cite{Ch-Gi-Pe}]\label{ImuHmuHp} Suppose that $1< p<\infty $ and let
$\mu $ be a positive Borel measure on $[0,1)$. Then:
\begin{itemize}
\item[(i)] The operator $I_\mu $ is well defined on $H^p$ if and only
if $\mu $ is a $1$-Carleson measure for $H^p$.
\item[(ii)] If $\mu $ is a $1$-Carleson measure for $H^p$, then the
operator $\mathcal H_\mu $ is also well defined on $H^p$ and,
furthermore,
$$\mathcal H_\mu (f)\,=\, I_\mu (f),\quad\text{for every $f\in
H^p$}.$$
\item[(iii)] The operator $I_\mu $ is a bounded operator from $H^p$
into itself if and only if $\mu $ is a Carleson measure.
\end{itemize}
\end{other}

Theorem\,\@\ref{ImuHmuH1} and Theorem\,\@\ref{ImuHmuHp} immediately
yield the following.

\begin{other}\label{HmuHp-prev} Let
$\mu $ be a positive Borel measure on $[0,1)$. \begin{itemize}
\item[(i)] If $\mu $ is a Carleson measure, then the
operator $\mathcal H_\mu $ is a bounded operator from $H^1$ into
itself if and only if  $\mu $ is a $1$-logarithmic $1$-Carleson
measure.
\item[(ii)]  If  $1<p<\infty $ and $\mu $ is a $1$-Carleson measure for $H^p$, then the
operator $\mathcal H_\mu $ is a bounded operator from $H^p$ into
itself if and only if $\mu $ is a Carleson measure.
\end{itemize}
\end{other}

Theorem\,\@\ref{HmuHp-prev} does not close completely the question
of characterizing the measures $\mu $ for which $\mathcal H_\mu $ is
a bounded operator from $H^p$ into itself. Indeed, in
Theorem\,\@\ref{HmuHp-prev} we only consider $1$-Carleson measures
for $H^p$. In principle, there could exist a measure $\mu $ which is
not a $1$-Carleson measures for $H^p$ but so that the operator
$\mathcal H_\mu $ is well defined and bounded on $H^p$. Our first
result in this paper asserts that this is not the case.

\begin{theorem}\label{HmuHp} Let
$\mu $ be a positive Borel measure on $[0,1)$. \begin{itemize}
\item[(i)] The operator $\mathcal H_\mu $ is a bounded operator from $H^1$ into
itself if and only if  $\mu $ is a $1$-logarithmic $1$-Carleson
measure.
\item[(ii)]
If  $1<p<\infty $ then the operator $\mathcal H_\mu $ is a bounded
operator from $H^p$ into itself if and only if $\mu $ is a Carleson
measure.
\end{itemize}
\end{theorem}
\par\smallskip We have the following result for $p=\infty $, a case
which was not considered in \cite{Ch-Gi-Pe}.
\begin{theorem}\label{HmuHinfty} Let
$\mu $ be a positive Borel measure on $[0,1)$. Then the following
conditions are equivalent.
\begin{itemize}
\item[(i)] $\int_{[0,1)}\frac{d\mu (t)}{1-t}\,<\,\infty $.
\item[(ii)] $\sum_{n=0}^\infty \mu_n\,<\,\infty $.
\item[(iii)] The operator $I_\mu $ is a bounded operator from
$H^\infty $ into itself.
\item[(iv)] The operator $\mathcal H_\mu $ is a bounded operator from
$H^\infty $ into itself.
\end{itemize}
\end{theorem}
\par\smallskip In the paper \cite{GM1} the authors have studied the
operators $\mathcal H_\mu $ acting on certain conformally invariant
spaces such as the Bloch space, BMOA, the analytic Besov spaces
$B^p$ ($1<p<\infty $), and the $Q_s$ spaces. Let us introduce
quickly these spaces.
\par It is well known that the set of all disc automorphisms ({\it
i.e.\/}, of all one-to-one analytic maps $f $ of $\mathbb D $ onto
itself), denoted $\Aut(\mathbb D) $, coincides with the set of all
M\"obius transformations of $\mathbb D $ onto itself: $ \Aut(\mathbb
D  )=\{ \lambda \varphi_a : \vert a \vert <1, \vert \lambda \vert =1
\} \,, $ where $\varphi_ a(z)=(a-z)/(1-\overline a z) $.
\par
A space $X $ of analytic functions in $\D $, defined via a semi-norm
$\rho $, is said to be {\it conformally invariant\/} or {\it
M\"obius invariant\/} if whenever $f\in X $, then also $f\circ
\varphi \in X $ for any $\varphi\in\operatorname{Aut}(\D) $ and,
moreover, $\rho (f\circ\varphi )\le C\rho (f) $ for some positive
constant $C $ and all $f\in X $. We mention \cite{AFP, DGV1,
Zhu-book} as references for M\"obius invariant spaces.
\par
The {\it Bloch space\/} $\mathcal B $ consists of all analytic
functions $f $ in $\mathbb D $ with bounded invariant derivative:
$$f\in \mathcal B \,\,\,\Leftrightarrow\,\,\,\rho _{\mathcal B}(f)\defeq  \sup\sb {z \in {\mathbb D
}}\,(1-|z|\sp2)\,|f\sp\prime(z)|<\infty \,. $$ A classical reference
for the Bloch space is \cite{ACP}; see also \cite{Zhu-book}. Rubel
and Timoney \cite{RT} proved that $\mathcal B$ is the biggest \lq\lq
natural\rq\rq\, conformally invariant space.
\par The space $BMOA $ consists of those functions $f $ in $H\sp1 $
whose boundary values have bounded mean oscillation on the unit
circle. Alternatively, $BMOA $ can be characterized in the following
way:
\par\smallskip
{\it If $f $ is an analytic function in $\mathbb D$, then $f\in BMOA
$ if and only if
$$\Vert f\Vert_{\star }
\,\defeq \,\sup_ {a\in \mathbb D }\Vert f\circ \varphi \sb
a-f(a)\Vert _ {H\sp 2}\,<\,\infty . $$} The seminorm $\Vert
\cdot\Vert_{\star }$ is conformally invariant. We mention
\cite{G:BMOA} as  a general reference for the space $BMOA$. Let us
recall that $$H^\infty \subsetneq BMOA \subsetneq
\bigcap_{0<p<\infty} H^p\quad \text{and}\,\, BMOA\subsetneq\mathcal
B.$$
\par\smallskip
If $0\le s<\infty$, we say that $f\in Q_s$ if $f$ is analytic in
$\D$ and \begin{equation*}\rho_{Q_s}(f)\,\defeq \,\left (
\sup_{a\in\D}\int_\D|f'(z)|^2g(z,a)^s\,dA(z)\right
)^{1/2}\,<\,\infty.\end{equation*} Here, $g(z,a)$ is the Green's
function in $\D$, given by
$g(z,a)=\log\left|\frac{1-\overline{a}z}{z-a}\right|$, while
$dA(z)=\frac{dx\,dy}{\pi }$ is the normalized area measure on
$\mathbb D$. All $Q_s$ spaces ($0\le s<\infty $) are conformally
invariant with respect to the semi-norm $\rho _{Q_s}$ (see e.\@g.,
\cite[p.\,\@1]{X2} or \cite[p.\,\@47]{DGV1}).
\par
These spaces were introduced by Aulaskari and Lappan
in~\cite{au-la94} while looking for new characterizations of Bloch
functions. They proved that for $s>1$, $Q_s$ is the Bloch space.
Using one of the many characterizations of the space $BMOA$ (see
\cite[Theorem\,\@6.\,\@2]{G:BMOA}) we have that $Q_1=BMOA$. In the
limit case $s=0$, $Q_s$ is the classical Dirichlet space $\mathcal
D$ of those analytic functions $f$ in $\D$ satisfying $\,
\int_\D|f'(z)|^2\,dA(z) <\infty$.
\par
Aulaskari, Xiao and Zhao proved in~\cite{au-xi-zh95} that
\[\mathcal D\subsetneq Q_{s_1}\subsetneq Q_{s_2}\subsetneq BMOA,\qquad0<s_1<s_2<1.\]
We mention \cite{X2} as an excellent reference for the theory of
$Q_s$-spaces.
\par\smallskip For $1<p<\infty $, the {\it analytic Besov space\/} $B\sp p $ is
defined as the set of all functions $f $ analytic in $\mathbb D $
such that
$$\rho_ p(f)=\left (\int_{\mathbb D}(1-\vert z\vert^2)^{p-2}\vert
f^\prime (z)\vert ^p\,dA(z)\right )^{1/p}\,<\,\infty .$$ All $B\sp p
$ spaces ($1<p<\infty $) are conformally invariant with respect to
the semi-norm $\rho _p$ (see \cite[p.\,\@112]{AFP} or
\cite[p.\,\@46]{DGV1}). We have that $\mathcal D\,=\,B^2$. A lot of
information on Besov spaces can be found in \cite{AFP, DGV1, DGV2,
HW1, Z, Zhu-book}. Let us recall that
$$B^p\,\subsetneq \,B^q\,\subsetneq BMOA,\quad
1\,<\,p\,<\,q\,<\infty  .$$
\par\smallskip
Among others, the following results have been proved in \cite{GM1}.
\begin{other}\label{ImuQs}
Let $\mu $ be a positive Borel measure on $[0,1)$.
\begin{itemize}\item[(i)] For any given $s>0$, the operator $I_\mu $ is
well defined in $Q_s$ if and only if
\begin{equation*}\label{log-finite}\int_{[0,1)}\log\frac{2}{1-t}d\mu
(t)<\infty .\end{equation*}
\item[(ii)] For any given $s>0$, the operator $I_\mu $ is a bounded
operator from $Q_s$ into $BMOA$ if and only if $\mu $ is a
$1$-logarithmic $1$-Carleson measure.
\item[(iii)] If $\mu $ is a
$1$-logarithmic $1$-Carleson measure then $\mathcal H_\mu (f)=I_\mu
(f)$, for all $f\in \mathcal B$.
\item[(iv)] If $\mu $ is a
$1$-logarithmic $1$-Carleson measure then $\mathcal H_\mu $ is a
bounded operator from $Q_s$ into $BMOA$ for any $s>0$.
\end{itemize}
\end{other}
\par It is natural to look for a characterization of those $\mu $
for which $I_\mu $ and/or $\mathcal H_\mu $ is a bounded operator
from $\mathcal B$ into itself or, more generally, from $Q_s$ into
itself for any $s>0$. We have the following result.
\begin{theorem}\label{ImuHmuQs} Let $\mu $ be a positive Borel measure on $[0,1)$.
Then the following conditions are equivalent.
\begin{itemize}
\item[(i)] The operator $I_\mu $ is bounded from $Q_s$ into itself
for some $s>0$.
\item[(ii)] The operator $I_\mu $ is bounded from $Q_s$ into itself
for all $s>0$.
\item[(iii)] The operator $\mathcal H_\mu $ is bounded from $Q_s$ into itself
for some $s>0$.
\item[(iv)] The operator $\mathcal H_\mu $ is bounded from $Q_s$ into itself
for all $s>0$.
\item[(v)] The measure $\mu $ is a
$1$-logarithmic $1$-Carleson measure.
\end{itemize}
\end{theorem}
\par\smallskip
In \cite{GM1} we also studied the operators $\mathcal H_\mu $ acting
on Besov spaces. Theorem\,\@3.\,\@8 of \cite{GM1} asserts that $\mu
$ being a $\gamma $-logarithmic $1$-Carleson measure for some
$\gamma
>1$ is a sufficient condition for the boundedness of $\mathcal H_\mu
$ from $B^p$ into itself, for any $p>1$. On the other hand,
Theorem\,\@3.\,\@7 of \cite{GM1} asserts that if $1<p<\infty $ and
the operator $\mathcal H_\mu $ is bounded from $B^p$ to itself then
$\mu $ is a $\gamma $-logarithmic $1$-Carleson measure for any
$\gamma <\,1-\frac{1}{p}$. We can improve this result as follows.
\begin{theorem}\label{Bp-log-nec-new} Suppose that $1<p<\infty $ and let
$\mu $ be a positive Borel measure on $[0,1)$ such that the operator
$\mathcal H_\mu $ is bounded from $B^p$ into itself. Then $\mu $ is
a $\left (1-\frac{1}{p}\right )$-logarithmic $1$-Carleson measure.
\end{theorem}
\par\smallskip The paper is organized as follows. The results
concerning Hardy spaces will be proved in
Section\,\@\ref{Hardy-spaces}; Section\,\@\ref{Mobius} will be
devoted to prove Theorem\,\@\ref{ImuHmuQs} and
Theorem\,\@\ref{Bp-log-nec-new}.
 We close this section
noticing that, as usual, we shall be using the convention that
$C=C(p, \alpha ,q,\beta , \dots )$ will denote a positive constant
which depends only upon the displayed parameters $p, \alpha , q,
\beta \dots $ (which sometimes will be omitted) but not  necessarily
the same at different occurrences. Moreover, for two real-valued
functions $E_1, E_2$ we write $E_1\lesssim E_2$, or $E_1\gtrsim
E_2$, if there exists a positive constant $C$ independent of the
arguments such that $E_1\leq C E_2$, respectively $E_1\ge C E_2$. If
we have $E_1\lesssim E_2$ and  $E_1\gtrsim E_2$ simultaneously then
we say that $E_1$ and $E_2$ are equivalent and we write $E_1\asymp
E_2$.
\section{The operator $\mathcal H_\mu $ acting on Hardy
spaces}\label{Hardy-spaces} This section is devoted to prove
Theorem\,\@\ref{HmuHp} and Theorem\,\@\ref{HmuHinfty}.\par\smallskip
\begin{Pf}{\,\em{Theorem \ref{HmuHp}\,\@(i).}}
Suppose that $\mathcal H_\mu $ is a bounded operator from $H^1$ into
itself. For $0<b<1$, set
$$f_b(z)=\frac{1-b^2}{(1-bz)^2},\quad z\in \mathbb D.$$
We have that $f_b\in H^1$ and $\Vert f_b\Vert_{H^1}=1$. Since
$\mathcal H_\mu $ is bounded on $H^1$, this implies that
\begin{equation}\label{H1fb}1\gtrsim \Vert \mathcal H_\mu (f_b)\Vert
_{H^1}.\end{equation} We also have,
$$f_b(z)\,=\,\sum_{k=0}^\infty a_{k,b}z^k,\quad\text{with
$a_{k,b}=(1-b^2)(k+1)b^k$.}$$ Since the $a_{k,b}$'s are positive, it
is clear that the sequence $\{ \sum_{k=0}^\infty \mu_{n+k}a_{k,b}\}
_{n=0}^\infty $ of the Taylor coefficients of $\mathcal H_\mu (f_b)$
is a decreasing sequence of non-negative real numbers. Using this,
Theorem\,\@1.\,\@1 of \cite{Pav-dec}, (\ref{H1fb}), and the
definition of the $a_{k,b}$'s, we obtain
\begin{align*} 1\,&\gtrsim \,\Vert \mathcal H_\mu (f_b)\Vert
_{H^1}\,\gtrsim \,\sum_{n=1}^\infty \frac{1}{n}\left
(\sum_{k=0}^\infty \mu_{n+k}a_{k,b}\right )\\ & =  \sum_{n=1}^\infty
\frac{1}{n}\left (\sum_{k=0}^\infty
a_{k,b}\int_{[0,1)}t^{n+k}\,d\,\mu (t)\right )\\
& \gtrsim \,(1-b^2)\sum_{n=1}^\infty \frac{1}{n}\left
(\sum_{k=1}^\infty kb^k\int_{[b,1)}t^{n+k}\,d\mu (t)\right )\\
& \gtrsim \,(1-b^2)\sum_{n=1}^\infty \frac{1}{n}\left
(\sum_{k=1}^\infty kb^{n+2k}\,\mu\left ([b,1)\right )\right )
\\ &
=\,(1-b^2)\mu \left ([b,1)\right )\sum_{n=1}^\infty
\frac{b^n}{n}\left (\sum_{k=1}^\infty kb^{2k}\right )
\\ & =
\,(1-b^2)\mu \left ([b,1)\right )\left (\log \frac{1}{1-b}\right
)\frac{b}{(1-b^2)^2}
\end{align*}
Then it follows that
$$\mu \left ([b,1)\right )\,=\,\og \left (
\frac{1-b}{\log\frac{1}{1-b}}\right ),\quad \text{as $b\to 1$}.$$
Hence, $\mu $ is a $1$-logarithmic $1$-Carleson measure.
\par The converse follows from Theorem\,\@\ref{HmuHp-prev}\,\@(i).
\end{Pf}
\par\smallskip
\begin{Pf}{\,\em{Theorem \ref{HmuHp}\,\@(ii).}} Suppose that
$1<p<\infty $ and that $\mu $ is a positive Borel measure on $[0,1)$
such that the operator $\mathcal H_\mu $ is a bounded operator from
$H^p$ into itself.
\par
For $0<b<1$, set
$$f_b(z)=\left (\frac{1-b^2}{(1-bz)^2}\right )^{1/p},\quad z\in \mathbb D.$$
We have that $f_b\in H^p$ and $\Vert f_b\Vert_{H^p}=1$. Since
$\mathcal H_\mu $ is bounded on $H^p$, this implies that
\begin{equation}\label{Hpfb}1\gtrsim \Vert \mathcal H_\mu (f_b)\Vert
_{H^p}.\end{equation} We also have,
$$f_b(z)\,=\,\sum_{k=0}^\infty a_{k,b}z^k,\quad\text{with
$a_{k,b}\thickapprox(1-b^2)^{1/p}k^{\frac{2}{p}-1}b^k$.}$$ Since the
$a_{k,b}$'s are positive, it is clear that the sequence $\{
\sum_{k=0}^\infty \mu_{n+k}a_{k,b}\} _{n=0}^\infty $ of the Taylor
coefficients of $\mathcal H_\mu (f_b)$ is a decreasing sequence of
non-negative real numbers. Using this, Theorem\,\@A of
\cite{Pav-dec}, (\ref{H1fb}), and the definition of the $a_{k,b}$'s,
we obtain
\begin{align*} 1\,&\gtrsim \,\Vert \mathcal H_\mu (f_b)\Vert
_{H^p}^p\,\gtrsim \,\sum_{n=1}^\infty n^{p-2}\left
(\sum_{k=0}^\infty \mu_{n+k}a_{k,b}\right )^p\\ & =
\sum_{n=1}^\infty n^{p-2}\left (\sum_{k=0}^\infty
a_{k,b}\int_{[0,1)}t^{n+k}\,d\mu (t)\right )^p\\
& \gtrsim \,(1-b^2)\sum_{n=1}^\infty n^{p-2}\left
(\sum_{k=1}^\infty k^{\frac{2}{p}-1}b^k\int_{[b,1)}t^{n+k}\,d\mu (t)\right )^p\\
& \gtrsim \,(1-b^2)\sum_{n=1}^\infty n^{p-2}\left (\sum_{k=1}^\infty
k^{\frac{2}{p}-1}b^{n+2k}\mu \left([b,1)\right )\right )^p\\ & =\,
(1-b^2)\mu\left ([b,1)\right )^p\sum_{n=1}^\infty n^{p-2}b^{np}\left
(\sum_{k=1}^\infty k^{\frac{2}{p}-1}b^{2k}\right )^p
\\ & \asymp \,(1-b^2)\mu\left ([b,1)\right )^p\frac{1}{(1-b)^2}\sum_{n=1}^\infty n^{p-2}b^{np}
\\ & \asymp \,\mu\left ([b,1)\right )^p\frac{1}{(1-b)^p},
\quad\text{as $b\to 1$}.\end{align*} Then it follows that $$\mu\left
([b,1)\right )\,=\,\og \left (1-b\right ),\quad\text{as $b\to 1$},$$
and, hence, $\mu $ is a Carleson measure.
\par The other implication follows from Theorem\,\@\ref{HmuHp-prev}\,\@(ii).
\end{Pf}
\par\medskip
\begin{Pf}{\,\em{Theorem \ref{HmuHinfty}.}}
The equivalence (i)\,$\Leftrightarrow$\,(ii) is clear because
$$\int_{[0,1)}\frac{d\mu (t)}{1-t}\,=\,\int_{[0,1)}\left
(\sum_{n=0}^\infty t^n\right )d\mu (t)\,=\,\sum_{n=0}^\infty
\int_{[0,1)}t^nd\mu (t)\,=\,\sum_{n=0}^\infty \mu_n.$$ The
implication (i)\,$\Rightarrow $\,(iii) is obvious.
\par (iii)\,$\Rightarrow $\,(i): Suppose (iii). Let $f$ be the
constant function $f(z)=1$, for all $z$. Then (iii) implies that
there exists a positive constant $C$ such that
$$\left \vert \int_{[0,1)}\frac{d\mu (t)}{1-tz}\right \vert \,\le
C,\quad z\in \mathbb D.$$ Taking $z=0$ in this inequality, (i)
follows.
\par (iii)\,$\Rightarrow $\,(iv): Suppose (iii). We have seen that then
(i) holds, and it is easy to see that (i) implies that $\mu $ is a
Carleson measure. Using part (ii) of Theorem\,\@\ref{ImuHmuH1}, it
follows that $\mathcal H_\mu $ is well defined in $H^\infty $ and
that $\mathcal H_\mu (f)=I_\mu (f)$ for all $f$ in $H^\infty $. Then
(iii) gives that $\mathcal H_\mu $ is bounded from $H^\infty $ into
itself.
\par (iv)\,$\Rightarrow $\,(iii): Suppose that (iv) is true and, as above,
let $f$ be the constant function $f(z)=1$, for all $z$. Then
$\mathcal H_\mu (f)\in H^\infty $. But $\mathcal H_\mu
(f)(z)=\sum_{n=0}^\infty \mu_nz^n$ and then it is clear that
$$\mathcal H_\mu (f)\in H^\infty \,\,\Leftrightarrow \,\,
\sum_{n=0}^\infty \mu_n<\infty .$$ Thus we have seen that
(iv)\,$\Rightarrow $\,(ii). Since (ii)\,$\Leftrightarrow$\,(iii),
this finishes the proof.
\end{Pf}

\section{The operator $\mathcal H_\mu $ acting on M\"{o}bius invariant spaces}\label{Mobius}
A basic ingredient in the proof of Theorem\,\@\ref{ImuHmuQs} will be
to have a characterization of the functions $f(z)=\sum_{n=0}^\infty
a_nz^n$ whose sequence of Taylor
 coefficients $\{ a_n\} _{n=0}^\infty $ is a decreasing sequence of
 nonnegative numbers which lie in the $Q_s$-spaces. This is quite
 simple for $s>1$ (recall that $Q_s=\mathcal B$ if $s>1$):\par Hwang and
 Lappan proved in \cite[Theorem\,\@1]{HL} that if $\{ a_n\} $ is a
 decreasing
sequence of nonnegative numbers then $f(z)=\sum_{n=0}^\infty a_nz^n$
is a Bloch function if and only if $a_n\,=\,\og \left
(\frac{1}{n}\right )$.
\par Fefferman gave a characterization of the analytic functions
having nonnegative Taylor coefficients which belong to $BMOA$,
proofs of this criterium can be found in \cite{Bon, G:BMOA, HoWa-Fe,
Sledd-Steg}. Characterizations of the analytic functions having
nonnegative Taylor coefficients which belong to $Q_s$ ($0<s<1$) were
obtained in \cite[Theorem\,\@1.\,\@2]{au-steg-xi96} and
\cite[Theorem\,\@2.\,\@3]{AGW}. Using the mentioned result in
\cite[Theorem\,\@1.\,\@2]{au-steg-xi96}, Xiao proved in
\cite[Corollary\,\@3.\,\@3.\,\@1, p.\,\@29]{X2} the following
result.
\begin{other}\label{dec-Qs}
Let $s\in (0, \infty )$ and let $f(z)=\sum_{n=0}^\infty a_nz^n$ with
$\{ a_n\} $ being a decreasing sequence of nonnegative numbers. Then
$f\in Q_s$ if and only if $a_n\,=\,\og \left (\frac{1}{n}\right )$.
\end{other}
\par Being based on Theorem\,\@1.\,\@2 of
\cite{au-steg-xi96}, Xiao's proof of this result is complicated. We
shall give next an alternative simpler proof. It will simply use the
validity of the result for the Bloch space and the simple fact that
the mean Lipschitz space $\Lambda ^2_{1/2}$ is contained in all the
$Q_s$ spaces ($0<s<\infty $) (see \cite[Remark\,\@4, p.\,\@427]{AGW}
or \cite[Theorem\,\@4.\,\@2.\,\@1.]{X2}).
\par We recall \cite[Chapter\,\@5]{D} that a function $f\in \hol (\mathbb D)$ belongs to
the mean Lipschitz space $\Lambda ^2_{1/2}$ if and only if
$$M_2(r,f^\prime )\,=\,\og \left (\frac{1}{(1-r)^{1/2}}\right ).$$
We have the following simple result for the space $\Lambda
^2_{1/2}$. \begin{lemma}\label{dec-lambda2} If $\{
a_n\}_{n=0}^\infty $ is a decreasing sequence of nonnegative numbers
and $f(z)=\sum_{n=0}^\infty a_nz^n$ ($z\in \mathbb D$), then $f\in
\Lambda^2_{1/2}$ if and only if\, $a_n\,=\,\og \left
(\frac{1}{n}\right ).$ \end{lemma} \begin{pf} If\, $a_n\,=\,\og
\left (\frac{1}{n}\right )$, then
$$M_2(r,f^\prime )^2\,=\,\sum_{n=1}^\infty n^2\vert a_n\vert^2
r^{2n-2}\,\lesssim \,\sum_{n=1}^\infty r^{2n-2}\,\lesssim
\,\frac{1}{1-r},$$ and, hence, $f\in \Lambda^2_{1/2}$.
\par Suppose now that $\{
a_n\}_{n=0}^\infty $ is a decreasing sequence of nonnegative numbers
and $f\in \Lambda^2_{1/2}$. Then, for all $n$
\begin{equation}\label{eee}\sum_{k=1}^nk^2a_k^2r^{2k-2}\,\le \sum_{k=1}^\infty
k^2a_k^2r^{2k-2}\,=\,M_2(r,f^\prime
)^2\,\lesssim\,\frac{1}{1-r}.\end{equation} Taking $r=1-\frac{1}{n}$
in (\ref{eee}), we obtain
\begin{equation}\label{eeee}\sum_{k=1}^nk^2a_k^2\,\lesssim\,n.\end{equation} Since
$\{ a_n\} $ is decreasing, using (\ref{eeee}) we have
$$a_n^2\sum_{k=1}^nk^2\,\lesssim\,\sum_{k=1}^nk^2a_k^2\,\lesssim\,n$$
and then it follows that $a_n\,=\,\og \left (\frac{1}{n}\right ).$
\end{pf}
\par\medskip Now Theorem\,\@\ref{dec-Qs} follows using the result of
Hwang and Lappan for the Bloch space, Lemma\,\@\ref{dec-lambda2},
and the fact that \begin{equation}\label{incl}\Lambda^2_{1/2}\subset
Q_s\subset \mathcal B,\quad\text{ for all $s$.}\end{equation}
\par\medskip Using (\ref{incl}), it is clear that Theorem\,\@\ref{ImuHmuQs}
follows from the following result. \begin{theorem}\label{th-X} Let
$\mu $ be a positive Borel measure on $[0,1)$ and let $X$ be a
Banach space of analytic functions in $\mathbb D$ with
$\Lambda^2_{1/2}\subset X\subset \mathcal B$. Then the following
conditions are equivalent.
\begin{itemize}
\item[(i)] The operator $I_\mu $ is well defined in $X$ and, furthermore, it is a bounded operator from $X$ into $\Lambda^2_{1/2}$.
\item[(ii)] The operator $\mathcal H_\mu $ is well defined in $X$ and, furthermore, it is a bounded operator from $X$ into $\Lambda^2_{1/2}$.
\item[(iii)] The measure $\mu $ is a
$1$-logarithmic $1$-Carleson measure.
\item[(iv)] $\int_{[0,1)}t^n\log\frac{1}{1-t}d\mu (t)\,=\,\og \left
(\frac{1}{n}\right )$.
\end{itemize}
\end{theorem}
\begin{pf}
According to Proposition\,\@2.\,\@5 of \cite{GM1}, $\mu $ is a
$1$-logarithmic $1$-Carleson measure if and only if the measure $\nu
$ defined by $d\nu (t)\,=\,\log\frac{1}{1-t}d\mu (t)$ is a Carleson
measure and, using Proposition\,\@1 of \cite{Ch-Gi-Pe}, this is
equivalent to (iv). Hence, we have shown that (iii)\,
$\Leftrightarrow $\, (iv).
\par
Set $F(z)\,=\,\log \frac{1}{1-z}$ ($z\in \mathbb D$). We have that
$F\in X$.
\par (i)\,$\Rightarrow $\,(iv): Suppose (i). Then
$$I_\mu (F)(z)\,=\,\int_{[0,1)}\frac{\log \frac{1}{1-t}}{1-tz}d\mu
(t)$$ is well defined for all $z\in \mathbb D$. Taking $z=0$, we see
that $\int_{[0,1)}\log \frac{1}{1-t}d\mu (t)<\infty $. Since $F\in
X$ we have also that $I_\mu (F)\in \Lambda ^{2}_{1/2}$, but
$$I_\mu (F)(z)\,=\,\int_{[0,1)}\frac{\log \frac{1}{1-t}}{1-tz}d\mu
(t)\,=\,\sum_{n=0}^\infty \left
(\int_{[0,1)}t^n\log\frac{1}{1-t}d\mu (t)\right )z^n.$$ Since the
sequence $\left \{ \int_{[0,1)}t^n\log\frac{1}{1-t}d\mu (t)\right \}
_{n=0}^\infty  $ is a decreasing sequence of nonnegative numbers,
using Lemma\,\@\ref{dec-lambda2} we see that (iv) holds.
\par (iv)\,$\Rightarrow $\,(i): Suppose (iv) and take $f\in X$.
Since $X\subset \mathcal B$, it is well known that $\vert f(z)\vert
\lesssim \log\frac{2}{1-\vert z\vert }$, see \cite[p.\,\@13]{ACP}.
This and (iv) give \begin{equation}\label{vertf}\int_{[0,1)}t^n\vert
f(t)\vert d\mu (t)\,=\,\og \left (\frac{1}{n}\right ).\end{equation}
Then it follows easily that $I_\mu (f)$ is well defined and that
$$I_\mu (f)(z)\,=\,\sum_{n=0}^\infty \left (\int_{[0,1)}t^nf(t)d\mu
(t)\right )z^n.$$ Now (\ref{vertf}) implies that
$\int_{[0,1)}t^nf(t)d\mu (t)\,=\,\og \left (\frac{1}{n}\right )$ and
then it follows that $I_\mu (f)\in \Lambda ^2_{1/2}$.
\par The implication (iv)\,$\Rightarrow $\,(ii) follows using
Theorem\,\@2.\,\@3 of \cite{GM1} and the already proved equivalences
(i)\,$\Leftrightarrow $\, (iii)\, $\Leftrightarrow$\, (iv).
\par It remains to prove that (ii)\,$\Rightarrow $\,(iv). Suppose
(ii) then $\mathcal H_\mu (F)\in \Lambda^{2}_{1/2}$. Now $$\mathcal
H_\mu (F)(z)\,=\,\sum_{n=0}^\infty \left (\sum_{k=1}^\infty\frac{
\mu_{n+k}}{k}\right )z^n.$$ Notice that the sequence $\{
\sum_{k=1}^\infty\frac{ \mu_{n+k}}{k}\} _{n=0}^\infty $ is a
decreasing sequence of nonnegative numbers. Then, using
Lemma\,\@\ref{dec-lambda2} and the fact that $\mathcal H_\mu (F)\in
\Lambda^{2}_{1/2}$, we deduce that
\begin{equation}\label{sum1n}\sum_{k=1}^\infty\frac{ \mu_{n+k}}{k}\,=\,\og \left
(\frac{1}{n}\right ).\end{equation} Now
$$\sum_{k=1}^\infty\frac{
\mu_{n+k}}{k}\,=\,\int_{[0,1)}\sum_{k=1}^\infty
\frac{t^{n+k}}{k}d\mu(t)\,=\,\int_{[0,1)}t^n\log\frac{1}{1-t}d\mu
(t).
$$ Then (iv) follows using (\ref{sum1n}).
\end{pf}
\par\smallskip
\begin{remark} It is clear that Theorem\,\@\ref{th-X} actually
implies the following result.
\begin{theorem}\label{refined} Let $\mu $ be a positive Borel
measure on $[0,1)$ and let $0<s_1, s_2<\infty $. Then following
conditions are equivalent. \begin{itemize}
\item[(i)] The operator $I_\mu $ is well defined in $Q_{s_1}$ and, furthermore, it is a bounded operator from $Q_{s_1}$ into $Q_{s_2}$.
\item[(ii)] The operator $\mathcal H_\mu $ is well defined in $Q_{s_1}$ and, furthermore, it is a bounded operator from $Q_{s_1}$ into $Q_{s_2}$.
\item[(iii)] The measure $\mu $ is a
$1$-logarithmic $1$-Carleson measure.
\end{itemize}
\end{theorem}
\end{remark}
\par\smallskip
\begin{Pf}{\,\em{Theorem \ref{Bp-log-nec-new}.}} Suppose that $1<p<\infty $ and let
$\mu $ be a positive Borel measure on $[0,1)$ such that the operator
$\mathcal H_\mu $ is bounded from $B^p$ into itself. For
$\frac{1}{2}<b<1$, set
$$g_b(z)\,=\,\left (\log \frac{1}{1-b^2}\right )^{-1/p}\log
\frac{1}{1-bz},\quad z\in \mathbb D.$$ We have,
$$g_b^\prime (z)\,=\,\left (\log \frac{1}{1-b^2}\right
)^{-1/p}\frac{b}{1-bz},\quad z\in \mathbb D$$ and then, using
Lemma\,\@3.\,\@10 of \cite{Zhu-book} with $t=p-2$ and $c=0$, we have
$$\int_{\mathbb D}(1-\vert z\vert ^2)^{p-2}\vert g_b^\prime (z)\vert
^p\,dA(z)\,\asymp \,\left (\log \frac{1}{1-b^2}\right
)^{-1}\int_{\mathbb D}\frac{(1-\vert z\vert ^2)^{p-2}}{\vert
1-bz\vert ^p}\,dA(z)\,\asymp\,1.$$ In other words, we have that
\begin{equation*}g_b\in B^p\quad\text{and}\quad \Vert g_b\Vert
_{B^p}\,\asymp\,1.\end{equation*} Since $\mathcal H_\mu $ is a
bounded operator from $B^p$ into itself, this implies that
\begin{equation}\label{gbBp}1\,\gtrsim \,\Vert \mathcal H_\mu
(g_b)\Vert _{B^p}^p.\end{equation} We have
$$g_b(z)=\sum_{k=0}^\infty a_{k,b}z^k,\quad\text{with}\quad
a_{k,b}\,=\,\left (\log \frac{1}{1-b^2}\right
)^{-1/p}\frac{b^k}{k}.$$ Since the $a_{k,b}$'s are positive it
follows that the sequence $\{ \sum_{k=0}^\infty \mu_{n+k}a_{k,b}\}
_{n=0}^\infty $ of the Taylor coefficients of $\mathcal H_\mu (g_b)$
is a decreasing sequence of non-negative real numbers. Using this,
\cite[Theorem\,\@3.\,\@10]{GM1}, and (\ref{gbBp}) we see that
\begin{align*}1\,\gtrsim &\,\Vert \mathcal H_\mu
(g_b)\Vert _{B^p}^p\,\gtrsim \,\sum_{n=1}^\infty n^{p-1}\left
(\sum_{k=1}^\infty \mu_{n+k}a_{k,b}\right )^p
\\ = & \,\left (\log
\frac{1}{1-b^2}\right )^{-1}\sum_{n=1}^\infty n^{p-1}\left
(\sum_{k=1}^\infty \frac{b^k}{k}\int_{[0,1)}t^{n+k}d\mu (t)\right
)^p
\\ \ge &\,\left (\log \frac{1}{1-b^2}\right
)^{-1}\sum_{n=1}^\infty n^{p-1}\left (\sum_{k=1}^\infty
\frac{b^k}{k}\int_{[b,1)}t^{n+k}d\mu (t)\right )^p
\\ \ge &\, \left (\log
\frac{1}{1-b^2}\right )^{-1}\sum_{n=1}^\infty n^{p-1}\left
(\sum_{k=1}^\infty \frac{b^{n+2k}}{k}\right )^p\mu \left
([b,1)\right )^p
\\ = &\, \left (\log
\frac{1}{1-b^2}\right )^{-1}\sum_{n=1}^\infty n^{p-1}b^{np}\left
(\sum_{k=1}^\infty \frac{b^{2k}}{k}\right )^p\mu \left ([b,1)\right
)^p
\\ = &\, \left (\log
\frac{1}{1-b^2}\right )^{p-1}\frac{1}{(1-b^p)^p}\mu\left
([b,1)\right )^p \\ \asymp &\, \left (\log \frac{1}{1-b^2}\right
)^{p-1}\frac{1}{(1-b)^p}\mu\left ([b,1)\right )^p.
\end{align*}
Then it follows that $\mu\left ([b,1)\right)\,\lesssim
\,\frac{1-b}{\left (\log\frac{1}{1-b}\right )^{1-\frac{1}{p}}}.$
This finishes the proof.
\end{Pf}

\par
\end{document}